# In the Search for the Infinite Servers Queue with Poisson Arrivals Busy Period Distribution Exponential Behavior

## Manuel Alberto M. Ferreira

Instituto Universitário de Lisboa
BRU-IUL, ISTAR-IUL
Av. das Forças Armadas, 1649-026
LISBOA, Portugal
Email: manuel.ferreira@iscte.pt

## José António Filipe[*]

Instituto Universitário de Lisboa
BRU-IUL, ISTAR-IUL
Av. das Forças Armadas, 1649-026
LISBOA, Portugal
Email: jose.filipe@iscte.pt
*Corresponding author

**Abstract**

This paper purpose is to investigate exponential behavior conditions for the M|G|∞ queue busy period length distribution. It is presented a general theoretical result that is the basis of this work. The complementary analysis rely on the M|G|∞ queue busy period length distribution moments computation. In M|G|∞ queue practical applications - in economic, management and business areas - the management of the effective number of servers is essential since the physical presence of infinite servers is not viable and so it is necessary to create that condition through an adequate management of the number of servers during the busy period.

José António Filipe: Assistant Professor at ISCTE-IUL (Instituto Universitário de Lisboa), has his Habilitation in Quantitative Methods, PhD in Quantitative Methods, Master in Management and Graduation in Economics. His current interests include, among others, Mathematics and Statistics, Multi Criteria Decision Making Methods, Chaos Theory, Game Theory, Stochastic Processes - Queues and Applied Probabilities, Bayesian Statistics - Application to Forensic Identification, Applications to Economics, Management, Business, Marketing, Finance and Social Problems in general.
https://ciencia.iscte-iul.pt/public/person/jcbf

**1 Introduction**

In this paper it is possible to see that queues theory makes available the use of a set of tools providing a practical operations management techniques package. This package is often used in very different areas of research to model and to explain practical cases of real life. It allows to model for example:

- cases of unemployment, health, or computers' logistic,
- staffing problems determination, business scheduling or companies' inventory levels,
- the improvement of customers satisfaction,
- telecommunications situations or supermarkets' realities, as it may be seen in this paper.

Indeed, real world queuing systems may be found in many concrete areas as it is the specific cases of commercial queuing systems (in which commercial companies serve external customers – hairdressers, garages, supermarkets, etc), transportation service systems (where trucks, for example, are the customers or servers – waiting to be loaded, for example), business-internal service systems (customers are internal to the organization – computer support system, for example), social service systems (waiting lists in a hospital in a country's health system, for example).

By considering and understanding queues, by learning and using them and knowing how to manage queues through models and equations, queues theory may contribute to improve companies' organization, namely customers management or organizations' internal processes development in order to give companies competitive advantages.

At the M|G|∞ queue system customers arrive following a Poisson process at rate $\lambda$, upon its arrival receive immediately a service with time length distribution function $G(.)$ and mean $\alpha$. The traffic intensity is $\rho = \lambda\alpha$. The M|G|∞ queue busy period length distribution is called *B*, the distribution function *B(t)* and the probability density function *b(t)*.

Note that this queue system main characteristic is that when a customer arrives it is immediately served. For that, it is not mandatory the physical existence of infinite servers. What interests is that when a customer arrives it finds immediately an available server. There are mainly two ways to do it:

- The customer is its own server, as it happens in a supermarket when customers are collecting the goods
- There is a bourse of servers that are made available when customers arrive. Of course it is necessary to dimension this bourse in order there is not customers lack. One example of this situation is the satellite processing information in war situations: when a message arrives it must be immediately processed.

These two examples show how important is this queue busy period length study because during that time it is fundamental to have either installation (in the first case) or servers (in the second case) available so that the system operates the best possible. Due to the exponential distribution interesting qualities, the "lack of memory" [1], the search for exponential behavior conditions for the M|G|∞ queue busy period length distribution is done through this work.

In the next section it is presented the main result on this search. Then from sections 3 till 7, the cases of some M|G|∞ queues, for particular service times, are considered. In section 8 the search of exponential behavior is performed through the moment's computation. This work ends with a conclusions section and a list of references.

## 2 The Main Result

The M|G|∞ queue busy period length Laplace-Stieltjes transform is, see [4],

$$\bar{B}(s) = 1 + \lambda^{-1}\left(s - \frac{1}{\int_0^\infty e^{-st - \lambda \int_0^t [1-G(v)]dv} dt}\right) \quad (2.1).$$

From (2.1) it is deduced

$$E[B] = \frac{e^\rho - 1}{\lambda} \quad (2.2),$$

for any the service time distribution.

**Proposition 2.1**

For service time distributions fulfilling $\lim_{\alpha \to \infty} G(t) = 0$, fixing $\lambda$, for $\alpha$ great enough $B$ is approximately exponential.

**Dem:** Inverting $\frac{1}{s}\bar{B}(s)$, with $\bar{B}(s)$ given by (2.1), it is obtained

$$B(t) = 1 - \lambda^{-1} \sum_{n=1}^\infty \left[\frac{e^{-\lambda \int_0^t [1-G(v)]dv} \lambda(1-G(t))}{1-e^{-\rho}}\right]^{*n} (1-e^{-\rho})^n, t \geq 0 \quad (2.3)$$

where $*$ is the convolution operator. Fixing $\lambda$, if $\lim_{\alpha \to \infty} G(t) = 0$, $1 - G(t) \cong 1$ for $\alpha$ great enough, and $\rho$ consequently great enough, then $\frac{e^{-\lambda \int_0^t [1-G(v)]dv} \lambda(1-G(t))}{1-e^{-\rho}} \cong e^{-\lambda t} \lambda$ and $B(t) \cong 1 - \lambda^{-1}\sum_{n=1}^\infty (\lambda e^{-\lambda t})^{*n}(1-e^{-\rho})^n$. This second member Laplace-Stieltjes transform is $\frac{1}{s}\frac{\lambda e^{-\rho}}{s+\lambda e^{-\rho}} + \frac{e^{-\rho}}{s+\lambda e^{-\rho}}$. And, after its inversion, it is got $1 - (1-e^{-\rho})e^{-\lambda e^{-\rho} t}$. So, it is concluded that $B(t) \cong 1 - (1-e^{-\rho})e^{-\lambda e^{-\rho} t}$. ∎

**Notes:**

---

[1]For instance, if the lifetime of a device is exponentially distributed, the probability that it goes on operational for a period of length $T_2$, after having been operating for a period of length $T_1$, is the same that if it had begun its operation:

$$P(T > T_1 + T_2 | T > T_1) = P(T > T_2).$$

That is: the device in $T_1$ has the same quality that in the beginning of the operation. From here the designation o "lack of memory "for this property. This is a reason for the importance of exponential distribution in reliability theory. It is a standard border between the situations $P(T > T_1 + T_2 | T > T_1) \leq P(T > T_2)$ and $P(T > T_1 + T_2 | T > T_1) \geq P(T > T_2)$ where the device incorporates negatively or positively, respectively, the effects of operation time (the memory).

- So, under Proposition 2.1 conditions, B is approximately exponentially distributed with mean $\frac{e^\rho}{\lambda}$, since for $\rho$ great enough $1 - e^{-\rho} \cong 1$. This is the main result of this study,
- This demonstration has the weakness of not giving any evaluation for the approximation error. This is an open problem.
- Also note that for $\alpha$ great enough it is negligible the difference between $\frac{e^\rho - 1}{\lambda}$ and $\frac{e^\rho}{\lambda}$,
- For a service time probability distribution supported only in the time interval $[a, b]$, $\lim_{\alpha \to \infty} G(t) = 0$ must be replaced by $\lim_{\alpha \to b} G(t) = 0$ in Proposition 2.1.

## 3 Special Service Times Distributions

Begin with

**Proposition 3.1**

For an M|G|∞ queue, if the service time distribution function belongs to the collection

$$G(t) = 1 - \frac{1}{\lambda} \frac{\left(1 - e^{-\rho}\right) e^{-\lambda t - \int_0^t \beta(u) du}}{\int_0^\infty e^{-\lambda w - \int_0^w \beta(u) du} dw - \left(1 - e^{-\rho}\right) \int_0^t e^{-\lambda w - \int_0^w \beta(u) du} dw},$$

$$t \geq 0, -\lambda \leq \frac{\int_0^t \beta(u) du}{t} \leq \frac{\lambda}{e^\rho - 1} \qquad (3.1)$$

the busy period length distribution function is

$$B(t) = \left(1 - (1 - G(0))\left(e^{-\lambda t - \int_0^t \beta(u) du} + \lambda \int_0^t e^{-\lambda w - \int_0^w \beta(u) du} dw\right)\right) *$$

$$* \sum_{n=0}^\infty \lambda^n (1 - G(0))^n \left(e^{-\lambda t - \int_0^t \beta(u) du}\right)^{*n}, -\lambda \leq \frac{\int_0^t \beta(u) du}{t} \leq \frac{\lambda}{e^\rho - 1} \quad (3.2). \blacksquare$$

**Notes:**

- The demonstration may be seen in [4],
- For $\frac{\int_0^t \beta(t) dt}{t} = -\lambda$, $G(t) = B(t) = 1, t \geq 0$ in (3.1) and (3.2), respectively,
- For $\frac{\int_0^t \beta(t) dt}{t} = \frac{\lambda}{e^\rho - 1}$, $B(t) = 1 - e^{-\frac{\lambda}{e^\rho - 1} t}, t \geq 0$, purely exponential, in (3.2)
- If $\beta(t) = \beta$ (constant)

$$G(t) = 1 - \frac{\left(1 - e^{-\rho}\right)(\lambda + \beta)}{\lambda e^{-\rho}\left(e^{(\lambda+\beta)t} - 1\right) + \lambda}, t \geq 0, -\lambda \leq \beta \leq \frac{\lambda}{e^\rho - 1} \qquad (3.3)$$

and

$$B^\beta(t) = 1 - \frac{\lambda + \beta}{\lambda}\left(1 - e^{-\rho}\right) e^{-e^{-\rho}(\lambda+\beta)t}, t \geq 0,$$

$$-\lambda \leq \beta \leq \frac{\lambda}{e^\rho - 1}, \qquad (3.4),$$

a mixture of a degenerate distribution at the origin and an exponential distribution,

- With $G(t)$ given by (3.1), $\int_0^\infty [1 - G(t)]dt =$

$$\int_0^\infty \left[ \frac{1}{\lambda} \frac{(1-e^{-\rho})e^{-\lambda t - \int_0^t \beta(u)du}}{\int_0^\infty e^{-\lambda w - \int_0^w \beta(u)du} dw - (1-e^{-\rho})\int_0^t e^{-\lambda w - \int_0^w \beta(u)du} dw} \right] dt = -\frac{1}{\lambda} \left[ \ln \left| \int_0^\infty e^{-\lambda w - \int_0^w \beta(u)du} dw - (1-e^{-\rho}) \int_0^t e^{-\lambda w - \int_0^w \beta(u)du} dw \right| \right]_0^\infty = \alpha,$$

such as it should happen, since $B$ is positively distributed,

- From (3.1), $\lim_{\alpha \to \infty} G(t) = 0, t \geq 0, -\lambda < \frac{\int_0^t \beta(u)du}{t} \leq \frac{\lambda}{e^\rho - 1}$. Then this service time distributions fulfill the Proposition 2.1 conditions. So, the distributions collection (3.2) have approximately exponential behavior for $\alpha$ great enough and $-\lambda < \frac{\int_0^t \beta(u)du}{t} \leq \frac{\lambda}{e^\rho - 1}$,

- For $\beta(t) = \beta$ (constant) it is easy to check directly that the approximately exponential behavior is assumed for $\alpha$ great enough.

To catch the meaning of $\alpha$ and $\rho$ great enough, it will be presented in the sequence the $B$ Coefficient of Variation, $\delta_1[B]$, Coefficient of Symmetry, $\delta_2[B]$ and Kurtosis, $\delta_3[B]$ computations for the systems $M|G_1|\infty$ (service time distribution given by (3.3) for $\beta = 0$). Note that for the $M|G_1|\infty$ queue the $E[B^n], n = 1,2, \ldots, n$ required to compute

Table 3.1: $M|G_1|\infty$

| Parameters $\rho$ | $\delta_1[B]$ | $\delta_2[B]$ | $\delta_3[B]$ |
|---|---|---|---|
| .5 | 2.0206405 | 9.5577742 | 15.983720 |
| 1 | 1.4710382 | 5.5867425 | 10.878212 |
| 10 | 1.0000454 | 4.0000000 | 9.0000000 |
| 20 | 1.0000000 | 4.0000000 | 9.0000000 |
| 50 | 1.0000000 | 4.0000000 | 9.0000000 |
| 100 | 1.0000000 | 4.0000000 | 9.0000000 |

these parameters are given by

$$E[B^n] = (1 - e^{-\rho}) \frac{n!}{(\lambda e^{-\rho})^n}, n = 1,2,\ldots \quad (3.5).$$

After $\rho = 10$ the $M|G_1|\infty$ busy period exponential behavior is evidenced (remember that for an exponential distribution $\delta_1[B] = 1, \delta_2[B] = 4$ and $\delta_3[B] = 9$).

## 4 Constant Service Times Distributions

Considering now constant (deterministic) service times, that is: the M|D|∞ queue, of course the constant distribution fulfils the conditions of Proposition 2.1. Continuing as in the former section to catch the exponential behavior, begin to note that (2.1) is equivalent, see [8], to $(\bar{B}(s)-1)(C(s)-1) = \lambda^{-1}sC(s)$ being $C(s) = \int_0^\infty e^{-st-\lambda\int_0^t[1-G(v)]dv}\lambda(1-G(t))dt$. Differentiating $n$ times, using Leibnitz's formula and making $s = 0$, it results

$$E[B^n] = (1)^{n+1}\left\{\frac{e^\rho}{\lambda}nC^{(n-1)}(0) - e^\rho\sum_{p=1}^{n-1}(-1)^{n-p}\binom{n}{p}E[B^{n-p}]C^{(p)}(0)\right\}, n = 1,2,... \quad (4.1)$$

with

$$C^{(n)}(0) = \int_0^\infty (-t)^n e^{-\lambda\int_0^t[1-G(v)]dv}\lambda(1-G(t))dt, n = 0,1,2,... \quad (4.2).$$

The expression (4.1) gives a recurrent method to compute $E[B^n], n = 0,1,2,...$ as a function of the $C^{(n)}(0), n = 0,1,2,...$. For the M|D|∞ system:

$$C^{(0)}(0) = 1 - e^{-\rho}$$
$$C^{(n)}(0) = -e^{-\rho}(-\alpha)^n - \frac{n}{\lambda}C^{n-1}(0), n = 1,2,... \quad (4.3)$$

and it is possible to compute any $E[B^n], n = 0,1,2,...$ exactly. So

Table 4.1: **M|D|∞**

| $\rho$ | $\delta_1[B]$ | $\delta_2[B]$ | $\delta_3[B]$ |
|---|---|---|---|
| .5 | .40655883 | 6.0360869 | 11.142336 |
| 1 | .56798436 | 4.5899937 | 9.6137084 |
| 10 | .99959129 | 4.0000000 | 9.0000000 |
| 20 | .99999999 | 4.0000000 | 9.0000000 |
| 50 | .99999999 | 4.0000000 | 9.0000000 |
| 100 | .99999999 | 4.0000000 | 9.0000000 |

and after $\rho = 10$ the M|D|∞ busy period exponential behavior is evidenced.

## 5 Exponential Service Times Distributions

Now, with exponential service times distribution, that is: for the M|M|∞ queue, begin to note that the exponential distribution fulfils the conditions of Proposition 2.1. To make a checking as for the M|D|∞ queue, it is not possible to obtain expressions as simple as (4.3) to the $C^{(n)}(0)$. It is mandatory to compute numerically integrals with infinite limits and so approximations must be done.

The results are:

Table 5.1: **M|M|∞**

| Parameters $\rho$ | $\delta_1[B]$ | $\delta_2[B]$ | $\delta_3[B]$ |
|---|---|---|---|
| .5 | 1.1109224 | 5.0972761 | 10.454678 |
| 1 | 1.1944614 | 5.4821324 | 10.923071 |
| 10 | 1.1227334 | 4.1511831 | 9.1617573 |
| 20 | 1.0544722 | 4.0326858 | 9.0337903 |
| 50 | 1.0206393 | 4.0049427 | 9.0550089 |
| 100 | 1.0101547 | 4.0012250 | 9.0012250 |

and only after $\rho = 20$ it can be said that those values are the ones of an exponential distribution.

## 6 Power Service Times Distribution

If the service distribution is a power function with parameter $c, c > 0$ $G(t) = \begin{cases} t^c & 0 \leq t < 1 \\ 1, & t \geq 1 \end{cases}$ and $\alpha = \frac{c}{c+1}$. So $\lim_{C \to \infty} G(t) = \begin{cases} 0, & 0 \leq t < 1 \\ 1, & t \geq 1 \end{cases}$ and $\lim_{C \to \infty} \alpha = 1$. Then it fulfils the conditions of Proposition 2.1, adapted in the last note. To check the busy period exponential behavior, in the usual form, for this system the values of $\delta_2[B]$ and $\delta_3[B]$ were computed for $\alpha = .25$, .5 and .8 making, in each case, $\rho$ assume values from .5 till 100. The results are:

Table 6.1: **M|P|∞**

| $\rho$ | $\alpha = .25$ | | $\alpha = .5$ | | $\alpha = .8$ | |
|---|---|---|---|---|---|---|
| | $\delta_2[B]$ | $\delta_3[B]$ | $\delta_2[B]$ | $\delta_3[B]$ | $\delta_2[B]$ | $\delta_3[B]$ |
| .5 | 3.0181197 | 9.5577742 | 1.5035507 | 5.9040102 | 3.8933428 | 9.3287992 |
| 1 | 4.4211164 | 9.1402097 | 2.7111584 | 7.4994861 | 3.9854257 | 9.0702715 |
| 1.5 | 5.3090021 | 10.433228 | 3.3711526 | 8.2784408 | 3.9749455 | 8.9969919 |
| 2 | 5.8206150 | 11.140255 | 3.7332541 | 8.6924656 | 3.9751952 | 8.9815770 |

| 2.5 | 6.0803833 | 11.489308 | 3.9322871 | 8.9173048 | 3.9809445 | 8.9828631 |
| 3 | 6.1786958 | 11.619970 | 4.0388433 | 9.0369125 | 3.9871351 | 8.9877124 |
| 6 | 5.7006232 | 11.020248 | 4.0969263 | 9.1024430 | 3.9996462 | 3.9996459 |
| 7 | 5.5034253 | 10.774653 | 4.0765395 | 9.0804332 | 3.9999342 | 8.9999341 |
| 8 | 5.3382992 | 10.570298 | 4.0596336 | 9.0623268 | 3.9999992 | 8.9999992 |
| 9 | 5.2037070 | 10.404722 | 4.0467687 | 9.0486468 | 4.0000086 | 9.0000086 |
| 10 | 5.0944599 | 10.271061 | 4.0372385 | 9.0385796 | 4.0000068 | 9.0000068 |
| 15 | 4.7702550 | 9.8790537 | 4.0152698 | 9.0156261 | 4.0000005 | 9.0000005 |
| 20 | 4.6102777 | 9.6888601 | 4.0082556 | 9.0083980 | 4.0000000 | 9.0000000 |
| 50 | 4.3045903 | 9.3338081 | 4.0012425 | 9.0012513 | 4.0000000 | 9.0000000 |
| 100 | 4.1715617 | 9.1842790 | 4.0003047 | 9.0003057 | 4.0000000 | 9.0000000 |

The analysis of the results shows a strong trend of $\delta_2[B]$ and $\delta_3[B]$, to 4 and 9, respectively, after $\rho=10$. This trend is faster the greatest is the value of $\alpha$.

## 7 Pareto Service Times Distribution

In this section only the exemplification method is used. Consider a Pareto distribution such that
$1 - G(t) = \begin{cases} 1, & t < k \\ \left(\frac{k}{t}\right)^3, & t \geq k \end{cases}$, $k>0$. Then $\alpha = \frac{3}{2}k$. The values calculated for $\delta_2[B]$ and $\delta_3[B]$ with $\lambda=1$ and, so, $\rho = \alpha$ are and show a strong trend from $\delta_2[B]$ and $\delta_3[B]$ to 4 and 9, respectively, after $\rho=20$. This is natural because, in this case, the convergence of $\alpha$ to infinite imposes the same behaviour to $k$. And so, for the above presented distribution function, it results $\lim_{\alpha \to \infty} G(t) = 0$.

Table 7.1: **M|Pa|∞**

| $\alpha = \rho$ | $\delta_2[B]$ | $\delta_3[B]$ |
|---|---|---|
| .5 | 1028.5443 | 1373.4466 |
| 1 | 1474.7159 | 1969.0197 |
| 10 | 38.879220 | 54.896896 |
| 20 | 4.0048588 | 9.0049233 |
| 50 | 4.0000000 | 9.0000000 |
| 100 | 4.0000000 | 9.0000000 |

But, considering now a Pareto distribution, such that $1 - G(t) = \begin{cases} 1, & t < .4 \\ \left(\frac{.4}{t}\right)^\theta, & t \geq .4 \end{cases}$, $\theta > 1$, so $\alpha = \frac{.4\theta}{\theta-1}$ and the values obtained for $\delta_2[B]$ and $\delta_3[B]$ in the same conditions as the previous case are and do not go against the hypothesis of the existence of a trend from $\delta_2[B]$ and $\delta_3[B]$ values to 4 and 9, respectively, although much slower than in the previous case. But, now, the convergence of $\alpha$ to infinite implies the convergence of $\theta$ to 1. So $\lim_{\alpha \to \infty} G(t) = \begin{cases} 0, & t < .4 \\ 1 - \frac{.4}{t} & t \geq .4 \end{cases}$ and it is not possible to guarantee at all that for $\alpha$ great enough $1 - G(t) \cong 1$.

Table 7.2: **M|Pa|∞**

| $\alpha = \rho$ | $\delta_2[B]$ | $\delta_3[B]$ |
|---|---|---|
| .5 | 10.993704 | 16.675733 |
| 1 | 6.8553306 | 12.010791 |
| 10 | 4.5112470 | 9.5724605 |
| 20 | 4.4832270 | 9.5397410 |
| 50 | 4.4669879 | 9.5208253 |
| 100 | 4.4616718 | 9.5146406 |

## 8 Looking for M|G|∞ Queue Busy Period Exponential Behavior through Moments Comparison

Also the busy period moments $E[B^n], n = 1,2,..8$ were computed for the M|$G_1$|∞, M|D|∞ and M|M|∞ queue systems. The results are presented below having been considered $\rho = .2, 1, 10, 20, 50, 100$ and $\lambda = 1$. For contrast effects it were also computed the same orders moments for the exponential distribution with mean $\frac{e^\rho - 1}{\lambda}$.

Table 8.1: $\rho = .5; \lambda = 1$

| $E[B^n]$ | M|$G_1$|∞, | M|D|∞ | M|M|∞ | Exponential Distribution with Mean $\frac{e^\rho - 1}{\lambda}$ |
|---|---|---|---|---|
| 1 | .64872127 | .64872127 | .64872127 | .64872127 |
| 2 | 2.1391211 | .49039984 | .94021749 | .84167857 |
| 3 | 10.580443 | .45345725 | 2.123908 | 1.6380444 |
| 4 | 69.776809 | .52362353 | 6.481435 | 4.2505369 |
| 5 | 575.2154 | .74561912 | 24.83009 | 13.787069 |
| 6 | 5690.1909 | 1.2729348 | 114.3113 | 53.663788 |
| 7 | 65670.772 | 2.5362864 | 614.2686 | 243.68989 |
| 8 | 866182.39 | 5.7760128 | 3773.0385 | 1264.6954 |

Table 8.2: $\rho = 1; \lambda = 1$

| $E[B^n]$ | M\|$G_1$\|∞, | M\|D\|∞ | M\|M\|∞ | Exponential Distribution with Mean $\frac{e^\rho-1}{\lambda}$ |
|---|---|---|---|---|
| 1 | 1.7182818 | 1.71828187 | 1.7182818 | 1.7182817 |
| 2 | 9.3415481 | 3.90498494 | 7.1649255 | 5.9049849 |
| 3 | 76.178885 | 11.974748 | 43.251592 | 30.439285 |
| 4 | 828.30271 | 48.000932 | 358.65020 | 209.21308 |
| 5 | 11257.801 | 240.00691 | 3702.6601 | 1797.4352 |
| 6 | 183611.26 | 1440.0037 | 45803.547 | 18531.001 |
| 7 | 3493750.0 | 10079.998 | 660802.68 | 222890.38 |
| 8 | 75975977. | 80639.996 | 10894769. | 3063907.9 |

Table 8.3: $\rho = 10$; $\lambda = 1$

| $E[B^n]$ | M\|$G_1$\|∞, | M\|D\|∞ | M\|M\|∞ | Exponential Distribution with Mean $\frac{e^\rho-1}{\lambda}$ |
|---|---|---|---|---|
| 1 | 22025.46 | 22025.466 | 22025.46 | 22025.46 |
| 2 | $9.7028634\times 10^8$ | $9.6984181\times 10^8$ | $1.0964476\times 10^9$ | $8.7024229\times 10^8$ |
| 3 | $6.4115936\times 10^{13}$ | $6.4057725\times 10^{13}$ | $8.1873951\times 10^{13}$ | $6.4110115\times 10^{13}$ |
| 4 | $5.6489899\times 10^{18}$ | $5.6412982\times 10^{18}$ | $8.1515907\times 10^{18}$ | $5.6482206\times 10^{18}$ |
| 5 | $6.2213642\times 10^{23}$ | $6.2100718\times 10^{23}$ | $1.0144929\times 10^{24}$ | $6.2202345\times 10^{23}$ |
| 6 | $8.2220799\times 10^{28}$ | $8.2034292\times 10^{28}$ | $1.5150846\times 10^{29}$ | $8.2202137\times 10^{28}$ |
| 7 | $1.2677235\times 10^{34}$ | $1.2642735\times 10^{34}$ | $2.6398031\times 10^{34}$ | $1.2673782\times 10^{34}$ |
| 8 | $2.2338775\times 10^{39}$ | $2.2267865\times 10^{39}$ | $5.2565179\times 10^{39}$ | $2.2331677\times 10^{39}$ |

Table 8.4: $\rho = 20$; $\lambda = 1$

| $E[B^n]$ | M\|$G_1$\|∞, | M\|D\|∞ | M\|M\|∞ | Exponential Distribution with Mean $\frac{e^\rho-1}{\lambda}$ |
|---|---|---|---|---|
| 1 | $4.8516519\times 10^8$ | $4.8516519\times 10^8$ | $4.8516519\times 10^8$ | $4.8516519\times 10^8$ |
| 2 | $4.7077053\times 10^{17}$ | $4.7077053\times 10^{17}$ | $4.97111287\times 10^{17}$ | $4.7077053\times 10^{17}$ |
| 3 | $6.8520443\times 10^{26}$ | $6.8520443\times 10^{26}$ | $7.6403133\times 10^{26}$ | $6.8520443\times 10^{26}$ |
| 4 | $1.3297494\times 10^{36}$ | $1.3297494\times 10^{36}$ | $1.5656919\times 10^{36}$ | $1.3297494\times 10^{36}$ |
| 5 | $3.2257405\times 10^{45}$ | $3.2257405\times 10^{45}$ | $4.0106193\times 10^{45}$ | $3.2257405\times 10^{45}$ |
| 6 | $9.3901022\times 10^{54}$ | $9.3901022\times 10^{54}$ | $1.2328148\times 10^{55}$ | $9.3901022\times 10^{54}$ |
| 7 | $3.1890255\times$ | $3.1890255\times$ | $4.4211069\times 10^{64}$ | $3.1890255\times 10^{64}$ |
| | $10^{64}$ | $10^{64}$ | | |
| 8 | $1.2377634\times 10^{74}$ | $1.2377634\times 10^{74}$ | $1.8119914\times 10^{74}$ | $1.2377634\times 10^{74}$ |

Table 8.5: $\rho = 50$; $\lambda = 1$

| $E[B^n]$ | M\|$G_1$\|∞, | M\|D\|∞ | M\|M\|∞ | Exponential Distribution with Mean $\frac{e^\rho - 1}{\lambda}$ |
|---|---|---|---|---|
| 1 | $5.1847055\times 10^{21}$ | $5.1847055\times 10^{21}$ | $5.1847055\times 10^{21}$ | $5.1847055\times 10^{21}$ |
| 2 | $5.3762343\times 10^{43}$ | $5.3762343\times 10^{43}$ | $5.4883410\times 10^{43}$ | $5.3762343\times 10^{43}$ |
| 3 | $8.3622575\times 10^{65}$ | $8.3622575\times 10^{65}$ | $7.8395261\times 10^{65}$ | $8.3622575\times 10^{65}$ |
| 4 | $1.7342337\times 10^{88}$ | $1.7342333\times 10^{88}$ | $1.5741896\times 10^{88}$ | $1.7342337\times 10^{88}$ |
| 5 | $4.4957455\times 10^{110}$ | $4.4957455\times 10^{110}$ | $3.9512479\times 10^{110}$ | $4.4974455\times 10^{110}$ |
| 6 | $1.3985470\times 10^{133}$ | $1.3985470\times 10^{133}$ | $1.19012551\times 10^{133}$ | $1.3985470\times 10^{133}$ |
| 7 | $5.0757381\times 10^{155}$ | $5.0757381\times 10^{155}$ | $4.1821348\times 10^{155}$ | $5.0757381\times 10^{155}$ |
| 8 | $2.1052966\times 10^{178}$ | $2.1052966\times 10^{178}$ | $1.6795590\times 10^{178}$ | $2.1052966\times 10^{178}$ |

Table 8.6: $\rho = 100$; $\lambda = 1$

| $E[B^n]$ | M\|$G_1$\|∞, | M\|D\|∞ | M\|M\|∞ | Exponential Distribution with Mean $\frac{e^\rho - 1}{\lambda}$ |
|---|---|---|---|---|
| 1 | $2.6881171\times 10^{43}$ | $2.6881171\times 10^{43}$ | $2.6881171\times 10^{43}$ | $2.6881171\times 10^{43}$ |
| 2 | $1.4451948\times 10^{87}$ | $1.4451948\times 10^{87}$ | $1.4599447\times 10^{87}$ | $1.4451948\times 10^{87}$ |
| 3 | $1.1654558\times 10^{131}$ | $1.1654558\times 10^{131}$ | $1.083083\times 10^{131}$ | $1.1654558\times 10^{131}$ |
| 4 | $1.2531527\times 10^{175}$ | $1.2531527\times 10^{175}$ | $1.1226720\times 10^{175}$ | $1.2531527\times 10^{175}$ |
| 5 | $1.6843107\times 10^{219}$ | $1.6843107\times 10^{219}$ | $1.4546350\times 10^{219}$ | $1.6843107\times 10^{219}$ |
| 6 | $2.7155746\times 10^{263}$ | $2.7155746\times 10^{263}$ | $2.2617075\times 10^{263}$ | $2.7165746\times 10^{263}$ |
| 7 | $4.1026610\times 10^{307}$ | $4.1026610\times 10^{307}$ | $4.1026610\times 10^{307}$ | $4.1026610\times 10^{307}$ |
| 8 | The program failed this calculation | The program failed this calculation | The program failed this calculation | The program failed this calculation |

The results evidence the trend to the exponential behavior as $\rho$ increases.

# 9 Conclusions

In the nowadays socio-economic context, the methodologies associated to queues theory are very interesting once they permit to model a set of cases and to get a formal interpretation of specific contexts and phenomena, allowing also to get a correct understanding of these situations. The case of the specific cases referred in this paper are symptomatic of the advantages of the application of queue model systems, taking benefits from the formal application of statistics and stochastic processes, for example, in the area of telecommunications or in supermarkets' studies as mentioned above in this study on queues theory.

In the present study, it is possible to see that the exponential distribution is very simple and quite useful from a practical point of view. It has been frequently considered in queuing systems study. Conditions under which $B$ is exponentially distributed or approximately exponentially distributed for the M|G|∞ queue were derived.

Many quantities of interest in queues are insensible. This means that they depend on the service time distribution only by its mean. Thus it is indifferent which service distribution is being considered. But using those given by (3.3), result quasi-exponential or exponential busy periods. And, for these service distributions, all distributions related to the busy period have simple forms and are related to the exponential distribution.

In section 2, for a large class of distributions under conditions of heavy- traffic, it was proved that $B$ is approximately exponential irrespectively of the service time distribution.

But, for instance, if the service distribution is a power function, as it was seen, such conditions must be adapted. However, for $\alpha$ near 1 and $\lambda$ and $\rho$ great enough, it is possible to guarantee that $B$ is approximately exponentially distributed.

Also if service distribution is a Pareto one, adaptations are needed, as it was shown in section 8, to identify conditions to guarantee the busy period exponential behaviour. And a situation was examined for which it was not possible to guarantee at all that for $\alpha$ great enough $1-G(t) \cong 1$.

Finally, in section 8, comparing directly the busy period length moments with those of an exponential distribution with the same mean it is evidenced the trend to exponential behavior as the traffic intensity increases.